\documentstyle[amsfonts]{article}
\title{On finitization of the Gordon identities}
\author{B. Feigin, S. Loktev}
\date{}
\newtheorem{proposition}{Proposition}
\newtheorem{lemma}{Lemma}
\newtheorem{theorem}{Theorem}
\newtheorem{conjecture}{Conjecture}
\newtheorem{corollary}{Corollary}
\newenvironment{proof}{\par \medskip{\it Proof.}\
}{\rule{7pt}{7pt}\par\medskip}
\newenvironment{lproof}{\par \medskip{\it Proof of the lemma.}\
}{\rule{5pt}{5pt}\par\medskip}
\newcommand{\bin}[2]{\left[{#1}\atop{#2}\right]}
\begin{document}
\maketitle

\renewcommand{\thefootnote}{}
\footnote{This work is partially supported by grants RFBR 99-01-01169,
INTAS-OPEN-97-1312}
\def \CC {{\mathbb C}}
\def \RR {{\mathbb R}}
\def \ZZ {{\mathbb Z}}

\def \SL {\widehat{SL_2}}
\def \sl {\widehat{sl_2}}
\def \NN {\widehat{\cal N}}
\def \nn {\widehat{n}}
\def \ch {{\rm ch}}

\begin{abstract}

In this paper we prove the identity that generalizes the Andrews--Gordon
identity. 
%But the bosonic part of our formula differs from the bosonic
%part of the Andrews--Gordon identity. 
Also we discuss the relation of our
formula to the geometry of affine flag varieties and to the geometry of
polyhedra. 

\end{abstract}

{\large\bf 1. Introduction.}

Let $\SL$ be the group corresponding to the affine Lie algebra $\sl$. This
group is a central extension of the group of matrices $2\times 2$ with
determinant one over the ring of functions on a circle.
By $\NN$ denote the subgroup of $\SL$ which consists of unipotent
upper--triangular matrices, by $\nn$ denote the Lie algebra of
$\NN$. In other words, the group $\NN$ consists of currents on the circle
with values in the
maximal nilpotent subgroup of $SL_2$. The algebra $\nn$ is abelian and it
can be identified with the space of functions on the circle. Let us choose
a basis $\{ e_i \}$, $i\in \ZZ$ in $\nn$, such that the element $e_m$   
corresponds to the function ${\rm e}^{2\pi i m \phi}$, where the
coordinate $\phi$ is the angle parameter on the circle.
Let $Gr$ be the affine grassmanian for the group $\SL$, i. e. the quotient
of $\SL$ by a maximal parabolic subgroup $P$. Suppose that $e_i$ belongs
to the Lie algebra of $P$ if $i\le 0$.
By $p\in Gr$ denote the image of $P$ in the variety $Gr$.

\def \MM {{\cal M}}

In the paper \cite{FS} the closure of the $\NN$--orbit of the point $p$ on
$Gr$ is studied. Denote this variety by $\MM$. The geometric definition of
the variety $\MM$ is rather complicated, but the description of $\MM$ in
terms of its coordinate ring seems more simple and clear to us.

Recall that there exists a line bundle $\xi$ on the variety $Gr$ such
that the space of sections of $\xi^{\otimes k}$ is a realization of
the irreducible vacuum representation $L_{0,k}$ at level $k$ of the
algebra $\sl$.
More precisely, $H^0(Gr,\xi^{\otimes k})\cong (L_{0,k})^*$. So the
coordinate ring of the variety $Gr$ is the direct sum of the spaces
$(L_{0,k})^*$ with the natural product
$(L_{0,k_1})^*\otimes(L_{0,k_2})^*\to (L_{0,k_1+k_2})^*$.

The variety $\MM$ is a subvariety of $Gr$, so its coordinate ring is a
quotient of the coordinate ring of $Gr$. Namely, let $v_k$ be the highest
vector of $L_{0,k}$ with respect to the subgroup $P$. The space $L_{0,k}$ 
contains the subspace $D(k) = U(\nn) v_k$. Then the coordinate ring of
$\MM$
is the direct sum of the spaces $D(k)^*$ with the induced product.

In integrable representations of $\sl$ at level $k$ we have the
relation $e(t)^{k+1} = 0$, when $e(t) = \sum e_i t^{-i}$. In \cite{FS} it
is shown that this relation leads to the following functional realization
of the spaces $D(k)^*$.

$$D(k)^* = \Sigma_0(k) \oplus \Sigma_1(k) \oplus \Sigma_2(k) \oplus
\dots,$$
where $\Sigma_m(k)$ is the space of symmetric polynomials in $m$ variables
$t_1, \dots t_m$ which are divisible by the product $t_1\dots t_m$ and
have a zero if $t_1 = \dots = t_{k+1}$.

In this notation the product on the direct sum of $D(k)^*$ is defined as
the set of maps $\Sigma_{m_1}(k_1) \otimes \Sigma_{m_2}(k_2) \to
\Sigma_{m_1+m_2}(k_1+k_2)$. The image of a pair of polynomials $P_1(t_1,
\dots t_{m_1})$ and $P_2(t_1, \dots, t_{m_2})$ under this map is the
symmetrization with respect to the variables $t_1, \dots, t_{m_1+m_2}$ of
their product:
\begin{equation}\label{umn}
(P_1\circ P_2)(t_1, \dots, t_{m_1+m_2}) = {\rm
Symm}(P_1(t_1, \dots, t_{m_1})P_2(t_{m_1+1}, \dots, t_{m_1+m_2})).
\end{equation}

Let $N$ be an integer. Let $\nn_N$ be the subalgebra of $\nn$ spanned by
$e_i$ for $i>N$, let $\NN_N$ be the corresponding subgroup.
We want to define the variety $\MM_N$ as the quotient of the variety $\MM$
by the group $\NN_N$. From the geometric point of view, to do it one
should remove all the non--stable orbits from the variety $\MM$ 
and factorize the resulting space by the action of the group. From the
algebraic point of view, the coordinate ring of $\MM_N$ is the subring of
$\NN_N$--invariants (or, just in other words, $\nn_N$--invariants) in the
coordinate ring of $\MM$.

Thus the coordinate ring of $\MM_N$ is the direct sum of spaces
$D_N(k)^*$ dual to the spaces $D_N(k) = \nn_N \setminus D(k)$ of
coinvariants of $D(k)$ with respect to the algebra $\nn_N$.  
The product on these spaces is defined in the similar way.
Now define the variety $\MM_N$ as the projective spectrum of the graded
algebra $\bigoplus D_N(k)^*$. In this way we obtain a line bundle
$\xi_N$ on $\MM_N$, such that $H_0(\MM_N, \xi_N^{\otimes k}) \cong
D_N(k)^*$.

On the variety $Gr$, there is an action of the two--dimensional torus $T$,
whose tangent space in the unit is spanned by the Cartan element
$h_0$ and the energy element $L_0$. It is easy to show that the action
of $T$ on the varieties  $\MM$ and $\MM_N$ is also defined and the action
of $T$ on $D(k)$ and $D_N(k)$ is defined by  $[h_0/2, e_i] = e_i$,
$[L_0, e_i] = i e_i$.

Let
$$d_N(k;q,z) = {\rm Tr}_{D_N(k)} q^{L_0} z^{h_0/2}$$
be the Hilbert polynomial of the graded space $D_N(k)$. There are two ways
to compute it, and it is shown in \cite{FS} that
they both succeed for the spaces $D(k)$.
The first way is purely algebraic. It is based on the
functional realization of $D_N(k)^*$ and the Gordon filtration on the
space of symmetric polynomials. This computation was done in \cite{FF},
and the result is written in the left hand side of (\ref{tozhd}). The
second way is an application of the holomorphic Lefschetz formula to the
bundle $\xi_N^{\otimes k}$. Let us discuss it in detail.

It can be shown that for even $N$ the variety $\MM_N$ is non--singular
(but for odd $N>3$ it contains exactly one singular point) and $H^i(\MM_N,   
\xi_N^{\otimes k}) = 0$ for $i>0$. In this case the fixed points of the
action of the torus $T$ on $\MM_N$ can be enumerated by two non--negative
integers $m$ and $n$, such that $m+n \le N/2$. The fixed point $x_{0,0}$
is the image of the point $p\in Gr$ in $\MM_N$. For $t\in \RR$ consider
the point  
$$ x_{m,n}(t) = \prod_{i=0}^{m-1} \exp(te_{2i+1}) \prod_{i=0}^{n-1}
\exp(te_{N-2i})\cdot x_{0,0}.$$
Then the fixed point $x_{m,n}$ is the limit $\lim\limits_{t \to   
\infty} x_{m,n}(t)$. Adding  the contributions of these points up we
obtain
the right hand side of (\ref{tozhd}). This expression is also written in
\cite{W}.

So for even $N$ we have the following formula.
{\small
\begin{equation}\label{tozhd}
\sum_{v_1, \dots, v_k \in \ZZ_+}
q^{\sum\limits_{i,j} v_i v_j\min(i,j)}
z^{\sum\limits_i iv_i} \prod_{i=1}^k \bin{(N+1)i -
2\sum\limits_j v_j\min(i,j) + v_i}{v_i}_q =
\end{equation}
}
{\scriptsize
$$
= \sum_{{m,n\in \ZZ_+}\atop{m+n \le N/2}} \frac{z^{k(m+n)} q^{k(m^2+(N+1)n
-
n^2)}}
{\prod\limits_{i=1}^m (1-q^i)(1-q^{-i-m+1}z^{-1})
\prod\limits_{i=2m+1}^{N-2n} (1-q^i z) \prod\limits_{i=1}^n (1-q^{-i})
(1-q^{-N+i+n-2}z^{-1})},$$
}
where $\bin{n}{m}_q$ is the $q$--binomial coefficient and the notation
$\frac1{1-q^\alpha z^\beta}$ should be considered as $\sum\limits_{i\ge 0}
q^{i\alpha} z^{i\beta}$ when $\alpha, \beta \ge 0$ and as
$\left(-\sum\limits_{i<0} q^{i\alpha} z^{i\beta}\right)$ when $\alpha,
\beta \le 0$.

%The right hand side of this formula can be obtained in a way not
%involving geometric description of $\MM_N$. 
%There is a similar situation with the Weyl formula. 

The right hand side of this formula resembles the alternating sum in the
Weyl formula.
Recall that the Weyl
formula for characters of integrable representations of a Kac--Moody
algebra can be interpreted as the Lefschetz formula applied to bundles on
the corresponding flag variety. The fixed points of the action of the
torus on the flag variety can be enumerated by the elements of the Weyl
group, so they correspond to the extremal vectors of an integrable
representation. And the contribution of a fixed point to the Lefschetz
formula
can be interpreted as the character of the representation "in a small
vicinity" around the corresponding extremal vector. Let us describe it for
representations of $\SL$.

Consider the integrable representation $L_{l,k-l}$ of $\SL$, let  
$v_{k,l}$ be its highest vector. Consider the extremal vector
$w(v_{k,l})$, when
$w$ is an element of the Weyl group of $\SL$. Denote the degrees of
this vector by $c_q(w;k,l)$ and $c_z(w;k,l)$. 
We call the limit
\begin{equation}\label{okr}
q^{c_q(w;k,l)} z^{c_z(w;k,l)} \lim_{r\to \infty}
\left(q^{-c_q(w;rk,rl)}
z^{-c_z(w;rk,rl)} \ch L_{rl,r(k-l)}\right)
\end{equation}
the character of the representation $L_{l,k-l}$ "in a small vicinity"
around the vector $w(v_{k,l})$. It is easy to show that it coincides with  
the image of the character of the Verma module $M_{l,k-l}$ under the
action of $w$. To obtain the contribution of the vector $w(v_{k,l})$ we
should
expand the formula (\ref{okr}) in positive powers of $q$. In this way we
obtain the character of a certain Verma module with the sign equal to
$(-1)^{l(w)}$. The Weyl formula states that the character of the
representation $L_{l,k-l}$ is the sum of such contributions.
 
In this paper we use this idea to calculate the characters of the spaces
$D_N(k)$. Namely, we decompose the character of $D_N(k)$ into the sum of
characters "in small vicinities" of some vectors.
In this way we
generalize the formula (\ref{tozhd}) and prove the generalization by
combinatorial methods. Indeed, the bosonic part of our formula is also
the Lefschetz formula applied to bundles on a certain variety.

\medskip
{\large\bf 2.1}
Consider the space
$$ D_N(k,l,r) = \CC [e_1, \dots, e_N]/ I_N(k,l,r),\quad N \ge 2,$$
where the ideal $I_N(k,l,r)$ is generated by the elements 
$$ e_1^{l+1}, \quad e_N^{r+1}, \quad \sum_{{1\le \alpha_1, \dots,
\alpha_{k+1}
\le N}\atop {\alpha_1+ \dots+ \alpha_{k+1} = i}} e_{\alpha_1}\dots
e_{\alpha_{k+1}},\ \mbox{when}\ k+1\le i \le N(k+1).$$

In other words, $I_N(k,l,r)$ is generated by the elements
$$ e_1^{l+1}, \quad e_N^{r+1}, \quad e(t)^{k+1},\ t\in \CC,\ \mbox{when}\
e(t) =
\sum_{i=1}^N e_i t^i.$$

In particular, the space $D_N(k)$ discussed in the introduction is
isomorphic to $D_N(k,k,k)$.

On the space $D_N(k,l,r)$ we introduce a bi--grading $(\deg_q,
\deg_z)$. To do it we are bi--grading the algebra $\CC[e_1, \dots, e_N]$
by
setting $\deg_q e_i = i$, $\deg_z e_i = 1$. As the generators of
$I_N(k,l,r)$ are homogeneous, this bi--grading on $\CC[e_1, \dots,
e_N]$ defines a bi--grading on
$D_N(k,l,r)$. By $d_N(k,l,r;q,z)$ denote the Hilbert polynomial of this
space:
$$d_N(k,l,r;q,z) = \sum_{i,j\ge 0} q^i z^j\dim (D_N(k,l,r))^{i,j}.$$

\medskip
{\large\bf 2.2}
The space dual to $D_N(k,l,r)$ has the following functional realization.

Let $\Sigma^N_m$ be the space of symmetric polynomials of the form
\begin{equation}\label{vid}
t_1 t_2\dots t_m f(t_1, t_2, \dots,
t_m)
\end{equation} 
with degree not exceeding $N$ in each of $m$ variables.
By $\Sigma^N_m(k,l,r)$ denote the subspace of $\Sigma^N_m$ which consists
of polynomials of the form (\ref{vid}), such that

\begin{enumerate}
\item $f(t_1, \dots, t_m) = 0$ if $t_1=\dots = t_{k+1}$.

\item $f(t_1, \dots, t_m) = 0$ if $t_1= \dots = t_{l+1} = 0$.

\item A polynomial $f(t, \dots, t, t_{r+2}, t_{r+3}, \dots, t_m)$ has the
degree less then $(r+1)(N-1)$ with respect to the variable $t$.
\end{enumerate}

\begin{proposition}
The space $D_N(k,l,r)$ is dual to the direct sum 
$$\Sigma^N_0(k,l,r) \oplus \Sigma^N_1(k,l,r) \oplus \Sigma^N_2(k,l,r)
\oplus \dots.$$
In these settings the subspace $(D_N(k,l,r))^{i,j}$ is dual to the
subspace of $\Sigma^N_j(k,l,r)$ which consists of homogeneous polynomials
of degree $i$.
\end{proposition}

\def \PP {{\cal P}}

The construction of the duality is the following. Let $S^m_N$ be the space
of
homogeneous polynomial of homogeneous degree $m$ in $N$ variables $e_1,
\dots, e_N$. Define the map $\PP: (S^m_N)^* \to \Sigma^N_m$ by setting for
$\phi \in (S^m_N)^*$

$$\PP(\phi)(t_1, \dots, t_m) = \phi(e(t_1)\cdot \dots\cdot
e(t_m)).$$

To complete the prove of the proposition it is enought to observe that
this map is an isomorphism and that the restriction 1) on the
subspace $\Sigma^N_m(k,l,r)$
corresponds to 
factorizing the dual space by $e(t)^{k+1}$, the restriction 2) corresponds
to factorizing by $e_1^{l+1}$ and the restriction 3) corresponds to
factorizing by $e_N^{r+1}$.

\def \alg {D_N^\vee}

\medskip
{\large\bf 2.3}
Note that we can define an associative commutative product on the space
$\alg = \Sigma^N_m(k,l,r)$. Namely, the
product of polynomials $P_1 \in \Sigma^N_{m_1}(k_1,l_1,r_1)$ and $P_2 \in
\Sigma^N_{m_2}(k_2,l_2,r_2)$ can be defined by the formula (\ref{umn}).
It is clear that $P_1\circ P_2 \in \Sigma^N_{m_1+m_2}(k_1+k_2, l_1+l_2,
r_1+r_2)$.

The structure of commutative algebra on the space $\alg$ implies the
following statement. 

\def \MMP {\MM_N'}

\begin{proposition}
There exists a variety $\MMP$ together with line bundles $\xi_K$,
$\xi_L$, $\xi_R$ and vector fields $L_q$, $L_z$ on $\MMP$, such that the
space $D_N(k,l,r)$ is dual to the space $H^0(\MMP,
\xi_K^{\otimes k}\otimes \xi_L^{\otimes l}\otimes \xi_R^{\otimes r})$ and
$d_N(k,l,r;q,z)$ is the trace of the action of $q^{L_q} z^{L_z}$ on this
space.
\end{proposition}

\def \MMPt {{\widetilde{\MM}}_N'}

The construction of the variety $\MMP$ is the following. Let $S$ be the 
spectrum of the algebra $\alg$. The algebra $\alg= \bigoplus D_N(k,l,r)^*$
is polygraded by integers $k$, $l$ and $r$, so there is an action of
three--dimensional complex torus $T_3$ on $S$. Let $\MMPt$ be the maximal
subvariety of $S$ where the torus $T_3$ acts freely. Set $\MMP = \MMPt
/T_3$. Then we obtain the bundles $\xi_K$, $\xi_L$, $\xi_R$ as the bundles
associated to the corresponding one--dimensional representations of the
torus $T_3$.
\medskip

{\large\bf 2.4}
Using the Gordon filtration on the spaces $\Sigma^N_m$ (see \cite{FF},
\cite{FS}) we obtain the following formula that generalize the left hand
side of (\ref{tozhd}).

\def \ll {{\cal L}}
\def \rr {{\cal R}}

\begin{theorem}
Let $Q$ be the $k\times k$ matrix with $Q_{i,j} = \min(i,j)$; let
$\ll$ and $\rr$ be the vectors in $\RR^k$ such that 
$$\ll_i = \left\{
\begin{array}{ll}
i-l & i \ge l\\
0   & i \le l
\end{array}\right. \qquad
\rr_i = \left\{
\begin{array}{ll}
i-r & i \ge r\\
0   & i \le r
\end{array}\right.
$$
Then $d_N(k,l,r;q,z) =$

$$ = \sum_{{v=(v_1, \dots, v_k)}\atop{v_i \in \ZZ_+}}
q^{v^tQv + \ll^t v}
z^{\sum iv_i} \prod_{i=1}^k \bin{(N+1)i - (2Qv+\ll+\rr-v)_i}{v_i}_q,$$
where the symbol ${}^t$ is just the transposition and 
$$\bin{n}{m}_q = \frac{\prod\limits_{i=1}^n (1-q^i)}{\prod\limits_{i=1}^m
(1-q^i) \prod\limits_{i=1}^{n-m} (1-q^i)}.$$
  
\end{theorem}

The proof of this formula repeats the proof of similar statements in
\cite{FF} and \cite{FS}.

\medskip
{\large\bf 2.5}
Next we construct a basis in the space $D_N(k,l,r)$.

\begin{theorem}
The monomials of the form $e_1^{a_1}\dots e_N^{a_N}$ such that
\begin{equation}\label{usl}
a_1 \le l, \quad a_N \le r, \quad
a_i + a_{i+1} \le k \ \mbox{for any}\ 1\le i < N,
\end{equation}
form a basis in $D_N(k,l,r)$.
\end{theorem}

\begin{proof}

In \cite{FF} this statement is proved for $r=k$ and, as $D_N(k,l,0)
\cong D_{N-1}(k,l,k)$, for $r=0$. 

First of all, from the fact that the monomials of the form $e_1^{a_1}\dots
e_N^{a_N}$ with the restrictions $a_i+ a_{i+1} \le k$ span the space
$D_N(k,k,k)$ it follows that the monomials with the restrictions
(\ref{usl})
span $D_N(k,l,r)$. So to complete the proof it is sufficient to compare
the dimension of the space $D_N(k,l,r)$ and the number of elements of the
presumed monomial basis.

\def \td {\tilde{D}}

By $\td_N(k,l,r)$ denote the set of monomials with the restrictions
(\ref{usl}). Also let $d_N(k,l,r) = \dim D_N(k,l,r)$ and 
let $\tilde{d}_N(k,l,r)$ be the number of elements in $\td_N(k,l,r)$.
We have
\begin{lemma}\label{l}\
\def \theenumi {\roman{enumi}}
\def\labelenumi {\theenumi)}

\begin{enumerate}
\item $\tilde{d}_N(k,l,r) = \tilde{d}_N(k,l,r-1) + 
\tilde{d}_{N-1}(k,l,k-r)$.

\item $d_N(k,l,r) \le d_N(k,l,r-1)+ d_{N-1}(k,l,k-r)$.
\end{enumerate}
\end{lemma}

\begin{lproof}

(i) Let $m = e_1^{a_1}\dots e_N^{a_N} \in \td_N(k,l,r)$. Then either
$a_N<r$ and $m \in \td_N(k,l,r-1)$, or $a_N=r$ and $m' = e_1^{a_1}\dots
e_{N-1}^{a_{N-1}}\in \td_{N-1}(k,l,k-r)$. And conversely, for any element
of $\td_N(k,l,r-1)$ or $\td_{N-1}(k,l,k-r)$ we can find the
corresponding element in $\td_N(k,l,r)$.
 
\def \coker {{\rm Coker}}
\def \im {{\rm Im}}

(ii) Introduce the map $\tilde{i}_N(r): \CC[e_1, \dots, e_{N-1}]\to
D_N(k,l,r)$ defined by $\tilde{i}_N(r) (P) = P\cdot e_N^r$. It is easy to
see that the image of the ideal $I_{N-1}(k,l,k-r)$ under this map is zero,
so $\tilde{i}_N(r)$ induces the map $i_N(r):
D_{N-1}(k,l,k-r) \to D_N(k,l,r)$. In particular we obtain $\dim \im
(i_N(r)) \le d_{N-1}(k,l,k-r)$.

It is also clear that the image of the ideal $I_N(k,l,r-1)$ under the
natural projection $D_N(k,l,r) \to \coker (i_N(r))$ is zero, so $\dim
\coker (i_N(r)) \le d_{N}(k,l,r-1)$. As $\dim D_N(k,l,r) = \dim \im
(i_N(r)) + \dim \coker (i_N(r))$, we have the statement of the lemma.

\end{lproof}

From the equality $d_N(k,l,r) = \tilde{d}_N(k,l,r)$ for $r=0,k$ and also
for $N=2$ it follows that the inequality in statement (ii) of
lemma~\ref{l} is in fact an equality. So we have the equality $d_N(k,l,r)
= \tilde{d}_N(k,l,r)$ for all $r$ and the statement of the theorem.

\end{proof}

\begin{corollary}\label{rec}
$d_N(k,l,r;q,z) = d_N(k,l,r-1;q,z)+ (q^Nz)^r d_{N-1}(k,l,k-r;q,z)$.
\end{corollary}

\medskip
{\large\bf 2.6}
Let us present another corollary from the theorem. Let $P_N(k,l,r)$ be the
polyhedron in $\RR^N = \{ (x_1, \dots, x_N)\}$ defined by the inequalities

$$ x_1\le l; \quad x_N \le r;\quad x_i+x_{i+1}\le k\ \mbox{for any}\
1\le i < N.$$

Consider the functionals  $\varphi_z = \sum x_i$, $\varphi_q = \sum i\cdot
x_i$ on $\RR^N$.

\begin{corollary}
The dimension of the space $D_N(k,l,r)$ coincides with the number of
integer points of the polyhedron $P_N(k,l,r)$, and 
$$d_N(k,l,r;q,z) = \sum_{{\bf x} \in \ZZ^N\cap P_N(k,l,r)}
q^{\varphi_q ({\bf x})} z^{\varphi_z({\bf x})}.$$
\end{corollary}

Following the ideas described in \cite{KP} we decompose the polynomial
$d_N(k,l,r;q,z)$ into the sum of contributions of the vertices of the
polyhedron $P_N(k,l,r)$. Namely, we
propose (see theorem \ref{stp}) an identity of the form
$$d_N(k,l,r;q,z) = \sum_{M \mbox{\ is a vertex of \ } P_N(k,l,r)} d_M.$$

Let $M$ be a vertex. In a small vicinity of $M$ our polyhedron looks like
a cone over a certain $N-1$--dimensional polyhedron. Denote this cone by
$C(M)$. In our case for the vertex $M$ we set
$$ d_M = \sum_{{\bf x} \in \ZZ^N\cap C(M)}
q^{\varphi_q ({\bf x})} z^{\varphi_z({\bf x})}.$$

According to \cite{KP} we say that a vertex $M$ is simple if $C(M)$ is a
cone over a simplex. In this case there exists a simple formula for the
contribution. It is easy to show that for generic $k,l,r$ (namely, for
$k>l>0$,
$k>r>0$, $l\ne r$, $l+r \ne k$) the simple vertices of $P_N(k,l,r)$ are 
$$M_{m,n} = \left( \underbrace{l,\ \ k-l,\ \ l,\ \ k-l,\dots}_m  \
,\underbrace{0,\dots, 0}_{N-m-n},
\ \underbrace{\dots, r,\ \ k-r,\ \ r}_n \right).$$

In this formula $n+m \le N$, and if $n+m = N$ then the sum of $m$th and
$(m+1)$th coordinates should not exceed $k$.  

\begin{proposition}
Suppose that $k$, $l$ and $r$ are generic. Then any integer point of
$C(M_{m,n})$ can be obtained from $M_{m,n}$ by adding a linear combination
with non--negative integer coefficients of the vectors $v^1, \dots, v^N$,
when
{\footnotesize
$$v^i = \left\{  
\begin{array}{ll}
\left( \underbrace{ 0,\dots, 0}_{m-i},\
\underbrace{ -1,\ 1,\dots, (-1)^{i-1},\ (-1)^i}_i,\
\underbrace{0,
\dots, 0}_{N-m}\right) & \quad i\le m,\\
\left( \underbrace{0, \dots, 0}_{i-1},\  1,\  \underbrace{0, \dots,
0}_{N-i}\right) & \quad m < i \le N-n,\\
\left( \underbrace{0, \dots, 0}_{N-n},\ \underbrace{(-1)^{N-i+1},\
(-1)^{N-i},\dots, 1,\ -1}_{N-i+1}, \underbrace{0,\dots,
0}_{n+i-N-1}\right) & \quad N-n < i
\end{array}\right.
$$
}
\end{proposition}

So
$$d_{M_{m,n}} = q^{\varphi_q(M_{m,n})}  z^{\varphi_z(M_{m,n})}
\prod_{i=1}^N
\frac1{1-q^{\varphi_q(v^i)}z^{\varphi_z(v^i)}}.
$$

The quantities $\varphi_q(M_{m,n})$ and $\varphi_z(M_{m,n})$ are linear
functions in $k$, $l$, $r$. The quantities $\varphi_q(v^i)$ and 
$\varphi_z(v^i)$ are independent on $k$, $l$, $r$ and equal to the
following numbers.

{\footnotesize
\begin{center}
\begin{tabular}{|l|c|c|c|c|c|}
\hline
& \multicolumn{2}{|c|}{$m \ge i$}  & $N-n \ge i$ &
\multicolumn{2}{|c|}{$i > N-n$} \\
\cline{2-3} \cline{5-6}
& $i$ is even & $i$ is odd & $i>m$ & $N-i$ is even & $N-i$ is odd \\
\hline
$\varphi_z(v^i)$ & 0 & -1 & 1 & -1 & 0 \\
\hline
$\varphi_q(v^i)$ & $i/2$ & $(i-1)/2-m$ & $i$ & $n-1-(3N-i)/2$ &
$-(N-i+1)/2$
\\
\hline
\end{tabular}
\end{center}
}

Thus $d_{M_{m,n}} =$
{\scriptsize

\begin{equation}\label{formula}
= \frac{q^{\varphi_q(M_{m,n})} z^{\varphi_z(M_{m,n})}}
{\prod\limits_{i=1}^{[m/2]}(1-q^i)
\prod\limits_{i=[m/2]+1}^{m}(1-q^{-i}z^{-1})
\prod\limits_{i=m+1}^{N-n}(1-q^i z)
\prod\limits_{i=1}^{[n/2]}(1-q^{-i})
\prod\limits_{i=[n/2]+1}^{n}(1-q^{i-N-1} z^{-1})
}
\end{equation}
}

To consider this formula as a series in positive powers of $q$ and $z$ let
us set the term $\frac1{1-q^\alpha z^\beta}$ equal to $\sum\limits_{i\ge
0} q^{i\alpha} z^{i\beta}$ if $\alpha, \beta \ge 0$ and equal to
$\left(-\sum\limits_{i<0} q^{i\alpha} z^{i\beta}\right)$ if
$\alpha, \beta \le 0$  (these are all the possible cases). Denote the
resulting series by $d^{m,n}_N(k,l,r;q,z)$. Formula (\ref{formula})
defines it for any integer $k$, $l$, $r$.

\medskip
{\large\bf 2.7}
It turns out that in our case it is enought to sum up the contributions
of the simple vertices ($M_{m,n}$), that is, the sum of the contributions
of
the other vertices is zero.
Let us set
$$d^e_N(k,l,r;q,z) = \sum_{{m+n<N \ \mbox{\tiny or}}\atop
{m+n = N,\, m \
\mbox{\tiny is even}}} d^{m,n}_N(k,l,r;q,z),$$
$$d^o_N(k,l,r;q,z) = \sum_{{m+n<N \ \mbox{\tiny or}}\atop {m+n = N,\, m \
\mbox{\tiny is odd}}} d^{m,n}_N(k,l,r;q,z).$$

\begin{theorem}\label{stp}\

\def \theenumi {\roman{enumi}}
\def\labelenumi {(\theenumi)}
\begin{enumerate}

\item If $N$ is odd and $l\le r$ then $d_N(k,l,r;q,z) = d^o_N(k,l,r;q,z)$,

\item If $N$ is odd and $l\ge r$ then $d_N(k,l,r;q,z) = d^e_N(k,l,r;q,z)$,

\item If $N$ is even and $l+r\le k$ then $d_N(k,l,r;q,z) =
d^o_N(k,l,r;q,z)$,

\item If $N$ is even and $l+r\ge k$ then $d_N(k,l,r;q,z) =
d^e_N(k,l,r;q,z)$.

\end{enumerate}  
\end{theorem}

\begin{proof}
First, prove some lemmas.

\begin{lemma}\label{l1}\
\def \theenumi {\roman{enumi}}
\def\labelenumi {(\theenumi)}
\begin{enumerate}

\item $d^e_N(k,l,r;q,z) = d^e_N(k,l,r-1;q,z) + (q^Nz)^r
d^e_{N-1}(k,l,k-r;q,z)$
\item $d^o_N(k,l,r;q,z) = d^o_N(k,l,r-1;q,z) + (q^Nz)^r
d^o_{N-1}(k,l,k-r;q,z)$

\end{enumerate}
\end{lemma}

\begin{lproof}   
It follows from the identity
$$d^{m,n}_N(k,l,r;q,z) = d^{m,n}_N(k,l,r-1;q,z) + (q^Nz)^r
d^{m,n-1}_{N-1}(k,l,k-r;q,z)$$
for $n>0$ and the identity $d^{m,0}_N(k,l,r;q,z) =
d^{m,0}_N(k,l,r-1;q,z)$.
\end{lproof}

\begin{lemma}\label{l2}
$d^e_N(k,l,-1;q,z) = 0$ for odd $N$, Á $d^o_N(k,l,-1;q,z) = 0$ for
even $N$.     
\end{lemma}

\begin{lproof}
It follows from the identity
$$d^{m,2n}_N(k,l,-1;q,z) = -d^{m,2n+1}_N(k,l,-1;q,z).$$
\end{lproof}

\begin{lemma}\label{l3}
If $l+r = k$ then $d^e_{2N}(k,l,r;q,z) = d^o_{2N}(k,l,r;q,z)$.
\end{lemma}   

\begin{lproof}

We want to prove that $\sum\limits_{m=0}^{2N} (-1)^m
d^{m,2N-m}_{2N}(k,a,k-a;q,z) =
0$ for any $k$ and $a$.
Let
\def \Pc {\Phi^{\vee}}
 $$\Phi(n) = \prod_{i=1}^n \frac1{1-q^i}, \quad \Psi(m,n) = \prod_{i=m}^n
 \frac1{1-q^{-i}z^{-1}}, \quad \Pc(n) = \prod_{i=1}^n \frac1{1-q^{-i}}.$$

Then
$$d^{2s,2N-2s}_{2N}(k,a,k-a;q,z) = P_{q,z} \cdot \Phi(s) \Psi(s+1,N+s)
 \Pc(N-s),$$
 $$d^{2s+1,2N-2s-1}_{2N}(k,a,k-a;q,z) = P_{q,z} \cdot \Phi(s)
\Psi(s+1,N+s+1) \Pc(N-s-1),$$
where $P_{q,z}$ is the common numerator of these fractions that doesn't
depend on $s$. It easy to prove by induction that for $s<N$ we have

{\small
 $$\sum_{m=0}^{2s} (-1)^m d^{m,2N-m}_{2N}(k,a,k-a;q,z) =
 \frac{P_{q,z}}{1-q^{-N}}\Phi(s)\Psi(s+1,N+s)\Pc(N-s-1),$$
 $$\sum_{m=0}^{2s+1} (-1)^m d^{m,2N-m}_{2N}(k,a,k-a;q,z) =
 -\frac{P_{q,z}}{1-q^N}\Phi(s)\Psi(s+2,N+s+1)\Pc(N-s-1).$$
}

Using the second formula for $s=N-1$ we obtain the statement of the lemma.
\end{lproof}

Now prove the theorem. For $N=2$ the proof is
straightforward. 
Here we already have two cases depending on $l$ and $r$.
Suppose that all the statements of the theorem are proved for $N-1$ and
prove them for $N$.

Let $N$ be odd. By lemmas ~\ref{l1} and \ref{l2} we have $d^e_N(k,l,0;q,z)
= d^e_{N-1}(k,l,k;q,z)$. As $D_N(k,l,0) \cong D_{N-1}(k,l,k)$, the
induction hypothesis implies the equality $d_N(k,l,0;q,z) =
d^e_N(k,l,0;q,z)$. Induction on $r$ based on lemma~\ref{l1} and
corollary~\ref{rec} shows that $d_N(k,l,r;q,z) =
d^e_N(k,l,r;q,z)$ if $r\le l$, which is statement (ii).

Statement (i) follows from (ii) by applying the transformation of the
space $D_N(k,l,r)$ mapping $e_i$ to $e_{N-i+1}$.

Now let $N$ be even. The proof of statement (iii) is similar to the proof
of (ii). For $r=k-l$ statement (iv) follows from statement (iii)
and lemma~\ref{l3}. Then statement (iv) for an arbitrary $r \ge k-l$
can be proved by similar induction on $r$.

\end{proof}

\medskip
{\large\bf 2.8}
To obtain the right hand side of (\ref{tozhd}) consider the particular
case $l=r=k$. In this case the vertices $M_{2m,2n}$, $M_{2m-1,2n}$,
$M_{2m,2n-1}$ and $M_{2m-1,2n-1}$ glue into one vertex (let us denote it
by $M'_{m,n}$) and, moreover, for odd $N$ all the vertices of the form
$M_{n,N-n}$ glue into one vertex (let us denote it by $M'$).

We set
$$d_{M'_{m,n}} = \sum\limits_{i,j =0,1}
d_N^{2m-i,2n-j}(k,k,k;q,z),$$
assuming that $d_N^{i,j}(k,k,k;q,z)=0$ for $i<0$ and for $j<0$.
Then the following statement implies the right hand side of (\ref{tozhd}). 

\begin{proposition}
{\footnotesize
$$d_{M'_{m,n}} = \frac{z^{k(m+n)} q^{k(m^2+(N+1)n -
n^2)}}
{\prod\limits_{i=1}^m (1-q^i)(1-q^{-i-m+1}z^{-1})
\prod\limits_{i=2m+1}^{N-2n} (1-q^i z) \prod\limits_{i=1}^n (1-q^{-i})
(1-q^{-N+i+n-2}z^{-1})}$$
}
\end{proposition}

Also let $d_{M'}= \sum\limits_{i+j=N,\ i\, \mbox {\scriptsize is even}}
d^{i,j}_N(k,k,k;q,z)$ for odd $N$. We propose
the following statement that describes the contribution of the singular
point on $\MM_N$.

\begin{conjecture} Let $N$ be odd. Then

$$d_{M'} = \frac{P_N(q,z)}
{\prod\limits_{i=1}^N (1-q^{-i}z^{-1})},$$
when $P_N(q,z)$ is a polynomial.
\end{conjecture}

\end{document}